\documentclass[a4paper,11pt]{article} 
\usepackage{amsmath}
\usepackage{amssymb}
\usepackage{amscd}
\usepackage[american]{babel}
\usepackage{color} \usepackage{graphics} 
\usepackage{graphicx}
\newtheorem{definition}{Definition}
\newtheorem{proposition}{Proposition}
\newtheorem{algorithm}{Algorithm}       
\newtheorem{theorem}{Theorem}
\newtheorem{corollary}{Corollary}

\newtheorem{example}{Example}
\newtheorem{remark}{Remark}

\title{Ribbon tableaux, ribbon rigged configurations and
  Hall-Littlewood functions at roots of unity}

\author{Francois Descouens\\ Universit\'e de Marne-la-Vall\'ee, France} 
\date{}
\begin{document}

\maketitle

\noindent \textbf{Abstract}: Hall-Littlewood functions indexed by
rectangular partitions, specialized at primitive roots of unity, can
be expressed as plethysms. We propose a combinatorial proof of this
formula using Schilling's bijection between ribbon tableaux and ribbon
rigged configurations.\\


\section{Introduction}
In \cite{LLT93,LLT94}, Lascoux, Leclerc and Thibon proved several
formulas for Hall-Littlewood functions $Q^{'}_\lambda(X;q)$ with the
parameter $q$ specialized at primitive roots of unity. This formula
implies a combinatorial interpretation of the plethysms
$l_k^{(j)}\lbrack h_\lambda \rbrack$ and $l_k^{(j)}\lbrack e_\lambda
\rbrack$, where $h_\lambda$ and $e_\lambda$ are respectively products
of complete and elementary symmetric functions, and $l_k^{(j)}$ the
Frobenius characteristics of representations of the symmetric group
$\mathfrak{S}_k$ induced by a transitive cyclic subgroup of
$\mathfrak{S}_k$. However, the combinatorial interpretation of the
plethysms of Schur functions $l_k^{(j)}\lbrack s_\lambda \rbrack$
would be far more interesting.  This question led the same authors to
introduce a new basis $H_\lambda^{(k)}(X;q)$ of symmetric functions,
depending on an integer $k \ge 1$ and a parameter $q$, which
interpolate between Schur functions, for $k=1$, and Hall-Littlewood
functions $Q^{'}_\lambda(X;q)$, for $k \ge l(\lambda)$. These
functions were conjectured to behave similarly under specialization at
root of unity, and to provide a combinatorial expression of the
expansion of the plethysm $l_k^{(j)}\lbrack s_\lambda \rbrack$ in the
Schur basis for suitable values of the parameters. This conjecture has
been proved only for the stable case $k=l(\lambda)$, which reduces to
the previous result on Hall-Littlewood functions, and $k=2$, which
gives the symmetric and antisymmetric squares $h_2[s_\lambda]$ and
$e_2[s_\lambda]$. The proof given in \cite{CarreLeclerc} relies upon
the study of diagonal classes of domino tableaux, i.e. sets of domino
tableaux having the same diagonals. Carr\'e and Leclerc proved that
the cospin polynomial of such a class has the form $(1+q)^aq^b$,
and from this obtained the specialization $H_{\lambda \cup \lambda}^{(2)}(X;-1)$.\\
In \cite{Schilling2}, Schilling defines ribbon rigged configurations
and gives a statistic preserving bijection between this kind of rigged
configurations and ribbon tableaux corresponding to a product of row
partitions.  The aim of this note is to show that the result on
Hall-Littlewood functions at roots of unity follows from an explicit
formula for the cospin polynomials of certain diagonal classes of
ribbon tableaux, which turn out to have a very simple characterization
through Schilling's bijection. \\\\
\textbf{Acknowledgment} All the computations on ribbon rigged
configurations, ribbon tableaux and Hall-Littlewood functions have
been implemented by using MuPAD-Combinat (see \cite{HT} for details on
this package and more especially \cite{FD} for implementation of
symmetric functions).


\section{Basic definitions}
%
Consider $\lambda=(\lambda_1,\ldots,\lambda_p)$ a partition². We denote by
$l(\lambda)$ its length $p$, $\vert \lambda \vert$ its weight $\sum_i
\lambda_i$, and $\lambda^{'}$ its conjugate partition. A $k$-ribbon is a
connected skew diagram of weight $k$ which does not contain a 2$\times$2
square.  The first (northwest) cell of a $k$-ribbon is called the
\textit{head} and the last one (southeast) the \textit{tail}. By removing
successively $k$-ribbons from $\lambda$ such that, at each step, the remaining
shape is still a partition, we obtain a partition $\lambda_{(k)}$ independant
of the removing procedure. This partition is called the $k$-core of $\lambda$.
The $k$-quotient is a sequence of $k$ partitions derived from $\lambda$ (in
\cite{JK} this sequence is computed using a system of abacus).  Denote by
$\mathcal{P}$ the set of all partitions, $\mathcal{P}_{(k)}$ the set of all
$k$-cores (i.e all partitions from which we cannot remove any $k$-ribbon) and
$\mathcal{P}^{(k)}$ the cartesian product
$\mathcal{P}\times\ldots\times\mathcal{P}$ of length $k$. The quotient
bijection $\Phi_k$ is defined by
\begin{equation}
\begin{array}{cccc}
\Phi_k: & \mathcal{P} & \longrightarrow & \mathcal{P}_{(k)} \times \mathcal{P}^{(k)} \\
        &  \lambda    &    \longmapsto  & (\lambda_{(k)}\ ,\ \lambda^{(k)}) \ .
\end{array}
\end{equation}
%
\begin{remark}\label{remarqueQuotient}
 {\rm  For any partition $\lambda$, we denote by $l(\lambda)\lambda$ the
  partition
  $(l(\lambda)\lambda_1,\ldots,l(\lambda)\lambda_{l(\lambda)})$. The
  $l(\lambda)$-core of $l(\lambda)\lambda$ is empty and the
  $l(\lambda)$-quotient of $l(\lambda)\lambda$ is the following product of single row
  partitions
  \begin{equation}
  (l(\lambda)\lambda)^{(l(\lambda))}=((\lambda_{l(\lambda)}), \ldots, (\lambda_1))\ .
  \end{equation}}
\end{remark}
\subsection{Preliminaries on $k$-ribbon tableaux}
\noindent A $k$-ribbon tableau of shape $\lambda$ and weight
$\mu$ is a tiling of the skew diagram $\lambda / \lambda_{(k)}$ by
labelled $k$-ribbons such that 
\begin{itemize}
\item[1-] the head of a ribbon labelled $i$ must not be on the right of a ribbon labelled $j > i$,
\item[2-] the tail of a ribbon labelled $i$ must not be on the top of
  a ribbon labelled $j \ge i$.
\end{itemize}
We denote by ${\rm Tab}_{\lambda,\mu}^{(k)}$ the set of all $k$-ribbon
tableaux of shape $\lambda$ and weight $\mu$, and by ${\rm
  Tab}_{\lambda}^{(k)}$ the set of all $k$-ribbon tableaux of shape
$\lambda$ and evaluation any composition of $\vert \lambda \vert/k$.
%
\begin{figure}[!h]
\begin{center}
\includegraphics[width=4cm]{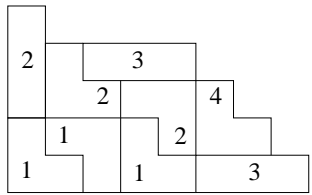}
\end{center}
\caption{A 3-ribbon tableau of shape (87651) and weight (3321)}
\end{figure}\\
%
The spin of a $k$-ribbon $R$ is defined by 
\begin{equation}
 {\rm sp}(R) = \frac{{\rm h}(R)-1}{2} \ ,
\end{equation}
where ${\rm h}(R)$ is the height of $R$. The spin of a $k$-ribbon
tableau $T$ of ${\rm Tab}_{\lambda,\mu}^{(k)}$ is the sum of the spins
of all its $k$-ribbons. The cospin of $T$ is the associated
co-statistic into ${\rm Tab}_{\lambda,\mu}^{(k)}$, i.e.
\begin{equation}
{\rm cosp}(T)=\max \left ({\rm sp}(T^{'}),\ T^{'} \in
  {\rm Tab}_{\lambda,\mu}^{(k)}\right ) - {\rm sp}(T)\ .
\end{equation}
We define the cospin polynomial $ \tilde G_{\lambda,\mu}^{(k)}(q)$ as
the generating polynomial of ${\rm Tab}_{\lambda,\mu}^{(k)}$ with respect to
the cospin statistic
$$
\tilde G_{\lambda,\mu}^{(k)}(q)=\sum_{T \in
 {\rm Tab}_{\lambda,\mu}^{(k)} } q^{{\rm cosp}(T)}~. $$
\begin{example} 
  The cospin polynomial for ${\rm Tab}_{(87651),(3321)}^{(3)}$ is
$$\tilde G_{(87651),(3321)}^{(3)}(q)= 3q^5+ 17q^4+33q^3+31q^2+18q +5~.$$
\end{example}
\subsection{$k$-tuples of semi-standard Young tableaux}
Let $T^{\circ}=(T^{\circ_1}, \ldots, T^{\circ_k})$ be a $k$-tuple of
semi-standard Young tableaux $T^{\circ_i}$ of shape
$\lambda^{\circ_i}$ and evaluation $\mu^{\circ_i}$. We call shape of
$T^\circ$ the sequence of partitions
$\lambda^{\circ}=(\lambda^{\circ_1},\ldots, \lambda^{\circ_k})$ and
weight of $T^\circ$, the composition $\mu=(\mu_1,\ldots,\mu_p)$
where $\mu_i=\sum_{j=1}^{k} \mu^{\circ_j}_i$.
%
\begin{definition}
Let $s$ be any cell in $T^\circ$
\begin{itemize}  
\item[-] \emph{pos($s$)} is the integer such that the cell $s$ is in the tableau
  $T^{\circ_{{\rm pos}(s)}}$,
\item[-] $T^\circ(s)$ is the label of the cell $s$,
\item[-] \emph{row($s$)} (resp. \emph{col($s$)}) represents the row
  (resp. the column) of $s$ in $T^{\circ_{{\rm pos}(s)}}$,
\item[-] \emph{diag($s$)} is the content of $s$, i.e. \emph{diag($s$)=col($s$)-row($s$)}.
\end{itemize}
\end{definition}
%
Denote by $\rm{Tab}_{\lambda^{\circ} ,\ \mu}$ the set of all $k$-tuples of
Young tableaux $T^\circ$ of shape $\lambda^\circ=(\lambda^{\circ_1},
\ldots, \lambda^{\circ_k})$ and weight $\mu$. In \cite{SW},
Stanton and White extend the bijection $\Phi_k$ to a correspondence
$\Psi_k$ between the set of $k$-ribbon tableaux of shape $\lambda$ and
weight $\mu$, and the set of $k$-tuples of semi-standard Young
tableaux of shape $\lambda^{(k)}$ and weight $\mu$.
%
\begin{example}
  The bijection $\Psi_3$ sends the $3$-ribbon tableau of Figure 1 to the 3-tuple of tableaux
$$
\left ( \ 
\begin{array}{|c|c|}
\hline 2 & 2 \\
\hline
\end{array}\quad , \quad 
\begin{array}{|c|c|c|}
\cline{1-2} 2 & 3 & \multicolumn{1}{c}{} \\
\hline 1 & 1 & 3 \\
\hline 
\end{array} \quad , \quad 
\begin{array}{|c|c|}
\hline 1 & 4 \\
\hline
\end{array}\ 
 \right ) \ .
$$
\end{example}
%
\begin{definition}[\cite{Schilling1}]Let $T^{\circ}=(T^{\circ_1}, \ldots, T^{\circ_k})$ 
  be a $k$-tuple of semi-standard Young tableaux and consider $s,t$
  two cells of $T^\circ$. The couple $(s,t)$ is an inversion in
  $T^\circ$ if the following conditions hold
\begin{itemize}
\item[1-] \emph{diag($s$)} = \emph{diag($t$)} and \emph{pos($s$)} $<$ \emph{pos($t$)} \\or \\
  \emph{diag($s$)} = \emph{diag($t$)}$ - 1$ and \emph{pos($s$)} $>$ \emph{pos($t$)},
\item[2-] \emph{row($s$)} $\le$ \emph{row($t$)},
\item[3-] $T(t) < T(s) < T(t^{\uparrow})$, where $t^{\uparrow}$ is the
  cell directly above $t$ and $T(t^{\uparrow})=\infty$ if
  $t^{\uparrow} \not \in \lambda^\circ$.
\end{itemize}
\end{definition}
The inversion statistic on $T^\circ$, denoted by Inv$(T^\circ)$, is
the number of couples in $T^\circ$ which form an inversion in
$T^\circ$.  This statistic permits to extend the correspondence
$\Psi_k$ to a bijection compatible with the inversion statistic and cospin
(\cite{Schilling1}), i.e.
\begin{equation}
\forall \ T \in {\rm Tab}_{\lambda,\mu}^{(k)}\ ,\ \text{Inv}\left
  (\Psi_k(T)\right )=\text{cosp}(T)\ .
\end{equation}
%
The inversion polynomial $\widetilde{I}_{\lambda^\circ,\mu}(q)$ is the
generating polynomial of ${\rm Tab}_{\lambda^{\circ},\mu}$ with
respect to the inversion statistic
\begin{equation}
\widetilde{I}_{\lambda^\circ,\mu}(q) = \sum_{T^\circ \in
  {\rm Tab}_{\lambda^\circ, \mu}} q^{{\rm Inv}(T^\circ)}\ .
\end{equation}
\begin{proposition}\label{I=G}
  The compatibility of the bijection $\Psi_k$ with the inversion
  statistic and cospin implies the following property
\begin{equation}\label{S=G}
 \widetilde{G}_{\lambda,\mu}(q) =\widetilde{I}_{\lambda^{(k)},\mu}(q)\ .
\end{equation}
\end{proposition}
In \cite{HHLRU}, Haglund, Haiman, Loehr, Remmel and Ulyanov use another notion
of inversion statistic. Their statistic coincides with Inv($T^\circ$) up to a
constant and gives a combinatorial interpretation of the powers of $q$
appearing in the decomposition of Macdonald polynomials
$\widetilde{H}_\lambda(X;q,t)$ on monomials.
\section{Specializations of Hall-Littlewood functions}
\subsection{Basic definitions}
Let $i,j$ be two nonnegative integers. The raising operator $R_{ij}$
acts on partitions by
\begin{equation}
\forall \ \lambda \in \mathcal{P},\ R_{ij}(\lambda)=(\lambda_1,\ldots,\lambda_i+1,\ldots,\lambda_j-1,\ldots,\lambda_p)\ .
\end{equation}
This action can be extended on elementary functions as follows
\begin{equation}
R_{ij}\cdot h_{\lambda} = h_{R_{ij}(\lambda)}\ .
\end{equation}
%
Hall-Littlewood functions $Q_{\lambda}^{'}(X;q)$ are
the symmetric functions defined by
\begin{equation}
Q_\lambda^{'}(X;q)=\prod_{i<j}(1-qR_{ij})^{-1}s_{\lambda}(X)\ .
\end{equation} 
We shall need the $\widetilde{Q}^{'}_\lambda(X;q)$ version of Hall-Littlewood
functions defined by
$$
\widetilde{Q}^{'}_\lambda(X;q)=q^{\eta(\lambda)}Q^{'}_\lambda(X;q^{-1})\ ,
$$
where $\eta(\lambda)=\sum_{i\ge 1}(i-1)\lambda_i$.
%
\begin{example}\label{exHL} 
  The expansion of the Hall-Littlewood function
  $\widetilde{Q}^{'}_{222}(X;q)$ on Schur basis is
\begin{eqnarray*}
\widetilde{Q}^{'}_{222}(X;q)&=&q^6s_{222}+(q^5+q^4)s_{321}+q^3s_{33}+q^3s_{411}+(q^4+q^3+q^2)s_{42}+(q^2+q)s_{51}+s_6\ .
\end{eqnarray*}
\end{example}
%
In \cite{LLT}, Lascoux, Leclerc and Thibon have shown that for any
partition $\lambda$, the Hall-Littlewood function
$\widetilde{Q}^{'}_\lambda(X;q)$ can be expressed in terms of
$l(\lambda)$-ribbon tableaux
$$
\widetilde{Q}_\lambda^{'}(X;q)~=\sum_{\mu\dashv \vert \lambda \vert}\ \sum_{T
  \in {\rm Tab^{(l(\lambda))}}_{l(\lambda)\lambda, \mu}}
q^{{\rm cosp}(T)}X^T = \sum_{\mu\dashv \vert \lambda \vert}
\widetilde{G}_{l(\lambda)\lambda,\mu}^{(l(\lambda))}(q)~m_{\mu}~,$$
where $X^T$ is
$x_1^{\mu_1}\ldots x_p^{\mu_p}$ with $\mu$ the weight of $T$.
\begin{proposition}\label{HLInversions}
  Hall-Litllewood functions $\widetilde{Q}_\lambda^{'}(X;q)$ can be
  expressed in terms of inversion polynomials using Proposition \ref{I=G}, 
  $$
  \widetilde{Q}_\lambda^{'}(X;q) = \sum_{\mu\dashv \vert \lambda \vert}
  \widetilde{I}_{(l(\lambda)\lambda)^{(l(\lambda))},\mu}(q)~m_{\mu} \ .
  $$
\end{proposition}
Using Remark \ref{remarqueQuotient}, this proposition means that the
Hall-Littlewood function $\widetilde{Q}^{'}_\lambda(X;q)$ can be
expressed in terms of $l(\lambda)$-tuples of Young tableaux with shape
the product of single row partitions
$((\lambda_{l(\lambda)}),\ldots,(\lambda_1))$.
%
\begin{example} 
  The expansion of the Hall-Littlewood function
  $\widetilde{Q}^{'}_{211}(X;q)$ on monomials is
  \begin{eqnarray*}
  \widetilde{Q}_{211}^{'}(X;q) &= & (t^3+3t^2+5t+3)m_{1111} + (t^3+t^2+2t)m_{22}+ \\
                               &  & (t^3+2t^2+3t+1)m_{211}+(t^3+t^2+t)m_{31} + t^3m_{4} \ .
\end{eqnarray*}
\end{example}
\subsection{Specialization at roots of unity}
Denote by $\Lambda_q$ the vector space of symmetric functions over the
field $\mathbb{C}(q)$. Let $k$ be a positive integer and $\lambda$ be a
partition. The plethysm of the powersum $p_\lambda$ by the powersum
$p_k$ is defined by
\begin{equation} 
p_k \circ p_\lambda = p_{k\lambda}\ .
\end{equation}
Since powersums $(p_\lambda)_{\lambda \in \mathcal{P}}$ form a
basis of the vector space $\Lambda_q$, the plethysm by a powersum
$p_k$ is defined on any symmetric function $f$.
\begin{theorem}[\cite{LLT93}] \label{LLT}
  Let $n,k$ be two positive integers and $\zeta$ be a primitive $k$-th
  root of unity. The specialization of the parameter $q$ at $\zeta$ in
  $Q^{'}_{n^k}(X;q)$ and $\widetilde{Q}^{'}_{n^k}(X;q)$ yields the
  following identities
\begin{equation}\label{HLRoots}
Q_{n^k}^{'}(X;\zeta) = (-1)^{(k-1)n}p_k\circ h_n(X)\quad {\rm and} \quad 
\widetilde{Q}_{n^k}^{'}(X;\zeta) = p_k\circ h_n(X)\ .
\end{equation}
\end{theorem}
\begin{example} 
  The Hall-Littlewood function $\widetilde{Q}^{'}_{222}(X;q)$ with $q$ specialized at
  $j=e^{\frac{2i\pi}{3}}$ is 
\begin{eqnarray}
\widetilde{Q}^{'}_{222}(X; j) & = &  s_{222}-s_{321} + s_{33} + s_{411} - s_{51} + s_6 \\
                  & = & \frac{1}{2}\left( p_{33} + p_{6} \right ) = p_3\circ \left (\frac{1}{2}\ p_{11} + \frac{1}{2}\ p_{2}\right )\\
                  & = & p_3\circ h_2(X)\ .
\end{eqnarray}
\end{example}
In \cite{LLT93}, Lascoux, Leclerc and Thibon have given an algebraic
proof of this theorem.
\section{Combinatorial proof of the specialization}
We give a combinatorial proof of (\ref{HLRoots}) for Hall-Littlewood functions
specialized at primitive roots of unity. The sketch of the proof is to
separate the set of $k$-tuples of single row Young tableaux into subsets
called diagonal classes. The specialization at primitive roots of unity of the
restriction of inversion polynomials on these classes are $0$ or $1$. In order
to have an explicit expression for these restricted polynomials, we translate
the problem into sets of ribbon rigged configurations according to
\cite{Schilling2}. These sets of configurations are interesting because the
image of diagonal classes can be easily characterized and inversion
polynomials have a nice expression in terms of fermionic formulas which behave
well at primitive roots of unity.
\noindent Due to Remark \ref{remarqueQuotient} and Proposition \ref{HLInversions}, we
only consider in the following, $k$-tuples of Young tableaux with shape
increasing sequences of row partitions $((\lambda_1^{\circ_1}),\ldots,
(\lambda_1^{\circ_p}))$, i.e.,
\begin{equation}\label{condition}
\lambda_1^{\circ_1} \le \ldots \le \lambda_1^{\circ_p}\ .
\end{equation}
%
\subsection{Diagonal classes of $k$-tuples of Young tableaux}

\indent Let $T^\circ=(T^{\circ_1},\ldots,T^{\circ_k})$ be a sequence of Young
tableaux with shape $\lambda^\circ=(\lambda^{\circ_1}, \ldots,
\lambda^{\circ_k})$ an increasing sequence of single row partitions.  Define
the maximal content $m$ of cells of $\lambda^\circ$ by
\begin{equation}
 m=\max\left (\lambda_1^{\circ_i}-1 , \ i={1\ldots k}\right) \ .
\end{equation}
For all $i\in \lbrace 0,\ldots, m\rbrace$, we
call $d_i$ the $i$-th diagonal of $T^\circ$ defined by
$$d_i=\lbrace T^\circ(s)\ \text{such that}\ \text{diag(s)}=i
\rbrace\ .$$
We call diagonal vector of $T^\circ$ the
vector $d_{T^\circ}=(d_0,\ldots,d_m)$. 
%
\begin{example}\label{exempleDiago}
The diagonal vector of the following 3-tuple of Young tableaux
$$\left ( \ 
\begin{array}{|c|c|}
\hline
1 & 4 \\
\hline
\end{array} \quad , \quad 
\begin{array}{|c|c|}
\hline
1 & 2 \\
\hline
\end{array} \quad , \quad
\begin{array}{|c|c|c|c|}
\hline
1 & 2 & 3 & 3 \\
\hline
\end{array}\ \right )  
$$
is given by $d_0=\{1,1,1\}$, $d_1=\{2,2,4\}$, $d_2=\{3\}$ and $d_3=\{3\}$.
\end{example}
%
\noindent Two $k$-tuples of Young tableaux $T^\circ$ and $T^{' \circ}$ in
${\rm Tab}_{\lambda^\circ,\mu}$ are equivalent if and only if, for all $i$ in
$\lbrace 0,\ldots,m\rbrace$, the $i$-th set in $d_{T^\circ}$ and $d_{T^{'
    \circ}}$ are the same. A diagonal class in ${\rm Tab}_{\lambda^\circ,\mu}$
is a set $D_{\lambda^\circ,\mu}(d)$ of all equivalent $k$-tuples of tableaux
with diagonal vector $d$.  We denoted by $\Delta_{\lambda^\circ, \mu}$ the set
of all diagonal vectors.  Thus, we can write the following decomposition
\begin{equation}\label{TabDiag}
{\rm Tab}_{\lambda^\circ, \mu} = 
\bigsqcup_{d\in \Delta_{\lambda^\circ,\mu}} D_{\lambda^\circ,\mu}(d) \ .
\end{equation}
%
The restriction of the inversion polynomial
$\widetilde{I}_{\lambda^\circ,\mu}(q)$ to a diagonal class
$D_{\lambda^\circ,\mu}(d)$ is defined by
\begin{equation}
\widetilde{I}_{\lambda^\circ,\mu}(q;d) = \sum_{T^\circ \in
  D_{\lambda^\circ,\mu}(d)}q^{{\rm Inv}(T^\circ)}\ .
\end{equation}
Hence, by (\ref{TabDiag}) the inversion polynomial $\widetilde{I}_{\lambda^\circ,\mu}(q)$
can be split into 
\begin{equation}\label{InvDiago}
\widetilde{I}_{\lambda^\circ,\mu}(q) = \sum_{d\in \Delta_{\lambda^\circ, \mu}}\widetilde{I}_{\lambda^\circ,\mu}(q;d) \ .
\end{equation}
%
\begin{example} 
  The diagonal vector of the 3-tuple of Young tableaux of Exemple
  \ref{exempleDiago} is $$d=(\{1,1,1\},\{2,2,4\},\{3\}, \{3\})\ .$$
  The diagonal class $D_{((2),(2),(4)), (3221)}(d)$ has 12 elements and the restriction of the
  inversion polynomial is
\begin{eqnarray*}
\widetilde{I}_{((2),(2),(4)), (3221)}(q;d) & = & q^5+3q^4+4q^3+3q^2+q \\
                                           & = & q(q+1)^2(q^2+q+1) \ .
\end{eqnarray*}
\end{example}
%
Since we are only interesting in Hall-Littlewood functions indexed by
partitions of the form $kn^k$, we will now restrict to the special case
of $k$-tuples of partitions $\lambda^\circ=\left ((n),\ldots,
  (n)\right )$ for some $n \ge 1$. We describe in the following
corollary the diagonal classes with only one element depending on the
weight $\mu$.
\begin{proposition}\label{one} 
  Let $n$ be a positive integer, $\lambda^\circ$ be the $k$-tuple of
  partitions $\left ((n),\ldots, (n)\right )$ and $\mu$ a partition of
  weight $kn$.
\begin{itemize}
\item[$\bullet$] If each part of $\mu$ is divisible by $k$
  \begin{itemize}
  \item[1-] there is a unique diagonal class
    $D_{\lambda^\circ,\mu}(d)$ with only one element
    $T^\circ=(T^{\circ_1},\ldots, T^{\circ_k})$,
  \item[2-] the $i$-th cell of each single row tableau $T^{\circ_j}$ is
    filled with the same value,
  \item[3-] $\emph{Inv}(T^\circ)=0$ and $\widetilde{I}_{\lambda^\circ, \mu}(q;d) = 1$.
\end{itemize}
\item[$\bullet$] if one part of $\mu$ is not divisible by $k$, there is no
  diagonal class with only one element.
\end{itemize}
\end{proposition}
\textit{Proof:} If each part of $\mu$ is divisible by $k$, we can construct a
diagonal vector $d$ such that for all $i$ in $\lbrace 0,\ldots, n-1\rbrace$,
all letters of $d_i$ are the same. And there is a unique way $T^\circ$ to fill
$\lambda^\circ$ according to $d$. Thus, the diagonal class $D_{\lambda^\circ,
  \mu}(d)$ has only one element.  Since each filling of $\lambda^{\circ_i}$
must be increasing, this
diagonal class is unique. This prove statement 1 and 2.\\
The inverse image of $T^\circ$ by Stanton-White map $\Psi_k$ is the $k$-ribbon
tableau of shape $nk^k$ containing $n$ blocks of size $k\times k$ of the form
\begin{center}
  \includegraphics*[width=3.5cm]{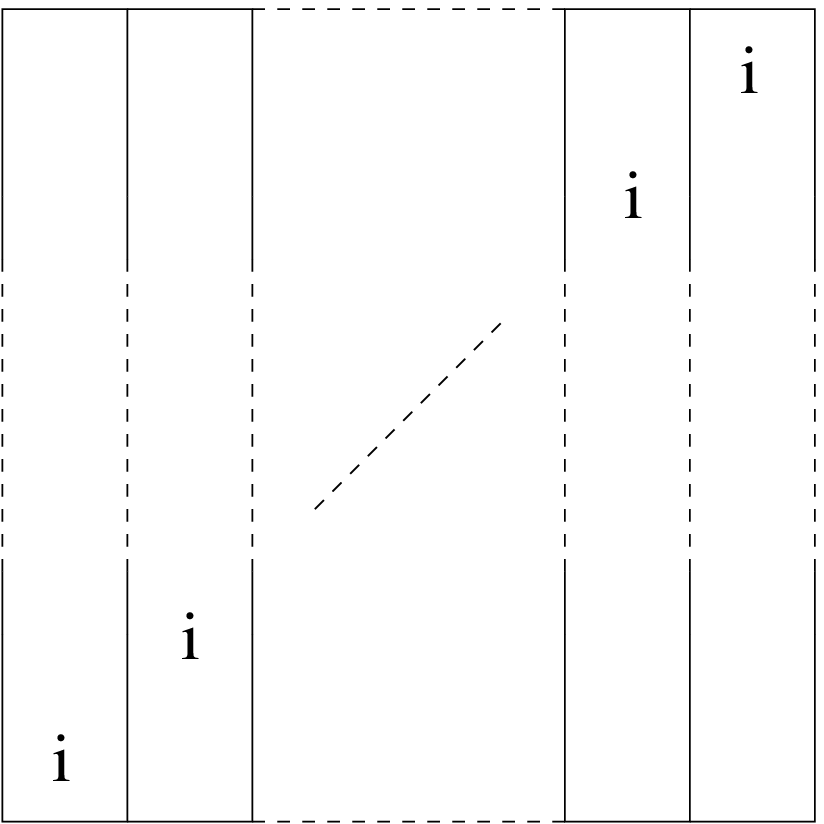}
\end{center}
Such a $k$-ribbon tableau has maximal spin in
$Tab^{(k)}_{\lambda,\mu}$. Thus, cospin of $T$ is equal to zero.
Since $\Psi_k$ is compatible with the inversion statistic and cospin, we
conclude that
$$
{\rm Inv}(T^\circ) = 0 \ .
$$
On the other hand, if one part of $\mu$ is not divisible by $k$,
this implies that for any filling $T^\circ$ of $\lambda^\circ$, there
exist two tableaux $T^{\circ_a}$ and $T^{\circ_b}$ in $T^\circ$ which
have two different values in one position. By transposition of these
two cells, we obtain another filling $T^{'\circ}$ which belongs to the
same diagonal class than $T^\circ$. Thus, for this kind of weight,
all diagonal
classes have more than 2 elements.
\hfill $\square$
\begin{proposition}\label{root} 
  Let $n$ be a positive integer, $\lambda^\circ$ be the $k$-uplet of
  partitions $((n),\ldots,(n))$ and $d$ be a diagonal vector such that
  $\#D_{\lambda^\circ, \mu}(d) \ge 2$. For any primitive $k$-th root of unity
  $\zeta$, the restriction of the inversion polynomial on
  $D_{\lambda^\circ,\mu}(d)$ satisfies the following specialization
\begin{equation}
  \widetilde{I}_{\lambda^\circ, \mu}(\zeta;d) = 0 \ . 
\end{equation}
\end{proposition}
In order to prove this proposition we need an explicit formula for the
polynomial $\widetilde{I}_{\lambda^\circ, \mu}(q;d)$. In \cite{Schilling2},
using a bijection between ribbon tableaux and ribbon rigged configurations,
Schilling gives a fermionic expression of the polynomial
$\widetilde{I}_{\lambda^\circ, \mu}(q)$ in terms of $q$-binomial coefficients.
This formula decomposes well on the image of diagonal classes on ribbon rigged
configurations.\\

%
\noindent One can define similarly diagonal classes on $k$-ribbon tableaux
throught the bijection $\Psi_k^{-1}$. The combinatorial interpretation of
cospin polynomials on diagonal classes is still an open problem for length of
ribbons $k>2$. For $k=2$ (case of domino tableaux), Carre and Leclerc have
found in \cite{CarreLeclerc} a combinatorial construction of these classes
using a notion of labyrinths and proved that cospin polynomials of diagonal
classes are of the form $q^a(1+q)^b$ with $a$ and $b$ two positive integers.
This paper give a solution for the stable case $k=l(\lambda)$ using ribbon
rigged configurations. We will come back to the proof of Proposition
\ref{root} in section \ref{preuveP4} after finishing the general proof of the
specialization of Hall-Littlewood functions.
%
\begin{corollary}\label{invpolroot}
  Let $\lambda^\circ$ be the $k$-tuple of single row partitions $((n), \ldots,
  (n))$, $\mu$ a partition of weight $nk$ and $\zeta$ be a $k$-th
  primitive root of unity. The specialization of the inversion
  polynomial at $q=\zeta$ yields
  \begin{itemize}
    \item[-] if all parts of $\mu$ are divisible by $k$
      $$
      \widetilde{I}_{\lambda^\circ, \mu}(\zeta) = 1 \ ,
      $$
    \item[-] if there exists a part of $\mu$ which is not divisible by $k$
      $$
      \widetilde{I}_{\lambda^\circ, \mu}(\zeta) = 0 \ .
      $$
  \end{itemize}
\end{corollary}
\textit{Proof:} By splitting the set $\Delta_{\lambda^\circ, \mu}$ with
respect to the cardinality of diagonal classes, (\ref{InvDiago}) can be
decomposed into
\begin{eqnarray}
\widetilde{I}_{\lambda^\circ, \mu}(q)& =& \sum_{d /
  \#D_{\lambda^\circ, \mu}(d)\ =\ 1} \widetilde{I}_{\lambda^\circ, \mu}(q; d) +
\sum_{d / \#D_{\lambda^\circ, \mu}(d)\ >\ 1 } \widetilde{I}_{\lambda^\circ,
  \mu}(q; d)\ .
\end{eqnarray}
By specializing $q$ at $\zeta$, the previous expression becomes
\begin{eqnarray}\label{eq}
\widetilde{I}_{\lambda^\circ, \mu}(\zeta)& =& \sum_{d /
  \#D_{\lambda^\circ, \mu}(d)\ =\ 1} \widetilde{I}_{\lambda^\circ, \mu}(\zeta; d) +
\sum_{d / \#D_{\lambda^\circ, \mu}(d)\ >\ 1 } \widetilde{I}_{\lambda^\circ,
  \mu}(\zeta; d)\ .
\end{eqnarray}
Using the result of Proposion \ref{one}, we conclude that the first
term gives 1 if all the parts of $\mu$ are divisible by $k$ and 0
otherwise. By Proposition \ref {root} the second term always gives 0,
which proves the corollary.  \hfill $\square$
\\\\
We are able to give a combinatorial proof of the specialization of
Hall-Littlewood functions given in Theorem \ref{LLT}. Let $n,k$ be two
positive integers, denote by $\Lambda_p^k$ the set of all partitions
of weight $p$ with all parts divisible by $k$. Using Corollary
\ref{invpolroot},
\begin{eqnarray}
\widetilde{Q}_{n^k}^{'}(X; \zeta) & = &\sum_{\mu\dashv nk} \widetilde{I}_{((n),\ldots, (n)), \mu}(\zeta)\ m_\mu \\
                      & = &\sum_{\mu\in \Lambda_{nk}^k}  m_\mu = \sum_{\mu\in \Lambda_{nk}^k}  p_k \circ m_{\mu/k} \ ,
\end{eqnarray}
where $\mu/k$ denote the partition
$(\frac{\mu_1}{k},\ldots,\frac{\mu_p}{k})$. The linearity of the
plethysm by $p_k$ implies
\begin{equation}
\widetilde{Q}_{n^k}^{'}(X; \zeta) = p_k\circ \left ( \sum_{\mu\in \Lambda_{nk}^k} m_{\mu/k} \right ) \ .
\end{equation}
By definition of complete functions, we conclude that
$$
\widetilde{Q}_{n^k}^{'}(X; \zeta) = p_k \circ h_n(X) \ .
$$
In rectangular case, the constant $\eta\left (n^k \right )$ is
equal to
$$\eta\left (n^k \right ) = n\sum_{i=1}^k(i-1)=\frac{nk(k-1)}{2}\ .$$
Hence, $\zeta^{\eta\left (n^k \right )} = (-1)^{(k-1)n}$ and
(\ref{HLRoots}) is proved.

%
\begin{remark} 
  {\rm In the case of $\lambda=k^{nk}$ which corresponds to
    $\lambda^{(k)}=((1^n),\ldots,(1^n))$, it exists
    a similar factorization formula than (\ref{HLRoots}) for Hall-Littlewood
    functions $Q^{'}_\lambda(X;q)$ and $\widetilde{Q}^{'}_\lambda(X;q)$
\begin{equation}
Q_{k^{nk}}^{'}(X;\zeta)= (-1)^{(k-1)n}p_k \circ e_n(X) \quad {\rm and} \quad
\widetilde{Q}_{k^{nk}}^{'}(X;\zeta)= p_k \circ e_n(X) \ .
\end{equation}
The proof is the same than for the product of single row partitions but using
ribbon rigged configurations defined for products of single column partitions
(\cite{Schilling2}).}
\end{remark}

\subsection{Ribbon rigged configurations and diagonal classes}
The aim of this section is to prove Proposition \ref{root} using
ribbon rigged configurations introduced by Schilling in
\cite{Schilling2}.\\\\
Let $n$ be a positive integer, we define the $q$-factorial $(q)_n$ by
\begin{equation}
(q)_n=\prod_{i=1}^{n}(1+q+\ldots+q^{i-1})\ .
\end{equation}
Let $a$ and $b$ be two positive integers. The $q$-binomial coefficient is
defined by
\begin{equation}
\left \lbrack
\begin{array}{c}
a+b \\ a , b
\end{array}
\right \rbrack = \frac{(q)_{a+b}}{(q)_a\ (q)_b} \ .
\end{equation}
\subsubsection{Definition of ribbon rigged configurations}
%
Let $\lambda, \mu$ be two partitions and
$\nu^\circ=\left(\nu^{\circ_1},\ldots,\nu^{\circ_p}\right )$ be a sequence of
partitions.  The sequence $\nu^\circ$ is a $(\lambda,\mu)$-configuration if
the following conditions hold 
\begin{itemize}
\item[1-] $\emptyset\subset \nu^{\circ_1} \subset \ldots \subset
  \nu^{\circ_{p-1}} \subset \nu^{\circ_p} = \widetilde{\mu}\ $ the conjugate partition of $\mu$,
\item[2-] $\forall\ 1 \le a < p, \quad \vert \nu^{\circ_a}\vert = \lambda_1 + \ldots + \lambda_a \ .$
\end{itemize}
We denote by $C_{\lambda, \mu}$ the set of all $(\lambda,
\mu)$-configurations. We associate to each element $\nu^\circ$ in
$C_{\lambda, \mu}$ a constant $\alpha(\nu^\circ)$ defined by
\begin{eqnarray}
\alpha(\nu^\circ) &=& \sum_{ {}_{1\le i \le \mu_1}^{1\le a \le p-1} }
\nu_{i+1}^{\circ_a}(\nu_i^{\circ_{a+1}}-\nu_i^{\circ_a})\ .
\end{eqnarray}
%
For any $(\lambda,\mu)$-configuration $\nu^\circ=(\nu^{\circ_1}, \ldots,
\nu^{\circ_p})$, the vacancy numbers $p_i^{(a)}$ and the constant
$m_i^{(a)}$ are defined for all $1\le a < p$ and $1 \le i \le \mu_1 $
by
\begin{equation}
p_i^{(a)} = \nu_i^{\circ_{a+1}}-\nu_i^{\circ_a} \quad {\rm and} \quad m_i^{(a)} = \nu_i^{\circ_a}-\nu_{i+1}^{\circ_a} \ .
\end{equation}
One can fill top cells of the $i$-th column of $\nu^{\circ_a}$ by a number
$j$ satisfying
\begin{equation}
\forall\ 1\le a < p,\ \forall\ 1 \le i \le \mu_1, \quad 0\le j \le p_i^{(a)} \ .
\end{equation}
Such a filling is called a rigging of $\nu^{\circ_a}$ and numbers $j$ are
called quantum numbers. In the special case of $j=p_i^{(a)}$, $j$ is called a singular quantum
number.\\
For each partition $\nu^{\circ_a}$ and height $i$, we can view a rigging as a
partition $J_i^{\circ_a}$ in a box of width $m_i^{(a)}$ and height
$p_i^{(a)}$. We denote by $J^{\circ_a}$ the $l(\nu^{\circ_a})$-tuple of
partitions $J^{\circ_a}=\left
  (J^{\circ_a}_1,\ldots,J^{\circ_a}_{l(\nu^{\circ_a})}\right )$ and by $J$ the
$(p-1)$-tuple $J=(J^{\circ_1},\ldots, J^{\circ_{p-1}})$. \\
Any riggings of $\nu^\circ$ by $J$ is called a rigged configuration $(\nu^\circ,J)$ of
shape $\nu^\circ$ and weight $J$.  The set of all rigged configurations $(\nu^\circ,
J)$ with $\nu^\circ \in C_{\lambda, \mu}$ is denoted by $RC_{\lambda,\mu}$.
\\\\
%
The graphical representation of a ribbon rigged configuration
$(\nu^\circ,J)$ is
\begin{itemize}
\item[1-] the filling of top cells of columns which are in the $i$-th row of
  $\nu^{\circ_a}$ with numbers of $J^{\circ_a}_i$ (riggings which differ only
  by reordering of quantum numbers corresponding to columns of the same height
  in a partition are identified),
\item[2-] the filling of the $i$-th row of the frontier of $\nu^{\circ_a}$
  with the vacancy numbers $p_i^{(a)}$.
\end{itemize}
%
\begin{example} The graphical representation of the ribbon rigged
  configuration in $RC_{(3111),(322)}$ given by
  $\nu^\circ=((21),(31),(32),(331))$ and $J=(((1),\emptyset), (\emptyset,(1)),
  (\emptyset,(1)))$ is
$$
\begin{array}{|c|c|c|}
\cline{1-1}
0 & \multicolumn{1}{c}{0}\\
\cline{1-2}
  & 1 &  \multicolumn{1}{c}{1}\\
\cline{1-2}
\end{array} \quad \quad 
\begin{array}{|c|c|c|c|}
\cline{1-1}
1 & \multicolumn{1}{c}{1}\\
\cline{1-3}
  & 0 & 0 & \multicolumn{1}{c}{1}\\
\cline{1-3}
\end{array} \quad \quad
\begin{array}{|c|c|c|c|}
\cline{1-2}
0 & 1&  \multicolumn{1}{c}{1}\\
\cline{1-3}
  &  & 0 & \multicolumn{1}{c}{0}\\
\cline{1-3}
\end{array} \quad \quad
\begin{array}{|c|c|c|}
\cline{1-1}
 & \multicolumn{2}{c}{}\\
\hline
 & &\\
\hline
  &  &  \\
\cline{1-3}
\end{array}\ .
$$
\end{example}
Denote by $RC_{\lambda, \mu}(\nu^\circ)$ the set of all rigged
configurations in $RC(\lambda, \mu)$ with shape $\nu^\circ$ in $C_{\lambda,\mu}$. Hence
\begin{equation}
  RC_{\lambda, \mu} =  \bigsqcup_{\nu\in C_{\lambda,\mu}} RC_{\lambda,
    \mu}(\nu^\circ)\ .
\end{equation}
There exists a cocharge statistic on ribbon rigged configurations defined by
\begin{eqnarray}
{\rm cc}(\nu^\circ, J)=\alpha(\nu^\circ) + \sum_{{}_{1\le i \le \mu_1}^{1\le a \le p-1}} \vert J_i^{\circ_a} \vert \ .
\end{eqnarray}
We denote by $\widetilde{S}_{\lambda, \mu}(q)$ the generating
polynomial of $RC_{\lambda,\mu}$ with respect to the cocharge, i.e.
\begin{eqnarray}\label{fermionic}
\widetilde{S}_{\mu, \delta}(q) & = & \sum_{(\nu,J) \in RC(\lambda, \mu)} q^{{\rm cc}(\nu^\circ,J)}
\end{eqnarray}
In \cite{Schilling2}, Schilling has given the following explicit expression
for $\widetilde{S}_{\lambda, \mu}(q)$,
\begin{eqnarray}
\widetilde{S}_{\lambda, \mu}(q) & = & \sum_{\nu^\circ \in C_{\lambda,\mu}} q^{\alpha(\nu^\circ)}
\prod_{ {}_{1\le i \le \mu_1}^{1\le a \le p-1} }\left\lbrack
{}_{\nu_i^{\circ_a}-\nu_{i+1}^{\circ_a},~
\nu_i^{\circ_{a+1}}-\nu_{i}^{\circ_a}}^{~~~~~~~\nu_i^{\circ_{a+1}}-\nu_{i+1}^{\circ_a} }
\right\rbrack\ .
\end{eqnarray}
This expression of $\widetilde{S}_{\lambda, \mu}(q)$ is called
fermionic formula. Denote by $\widetilde{S}_{\lambda, \mu}(q; \nu^\circ)$ the
cocharge polynomial $\widetilde{S}_{\lambda, \mu}(q)$ restricted to the
subset $RC_{\lambda, \mu}(\nu^\circ)$, i.e.
\begin{equation}\label{FermionicDiago}
\widetilde{S}_{\lambda, \mu}(q; \nu^\circ) =  q^{\alpha(\nu^\circ)}
\prod_{ {}_{1\le i \le \mu_1}^{1\le a \le p-1} }\left\lbrack
{}_{\nu_i^{\circ_a}-\nu_{i+1}^{\circ_a},~
\nu_i^{\circ_{a+1}}-\nu_{i}^{\circ_a}}^{~~~~~~~\nu_i^{\circ_{a+1}}-\nu_{i+1}^{\circ_a} }
\right\rbrack \ .
\end{equation}
\subsubsection{Bijection with $k$-tuples of tableaux}
\begin{theorem}[\cite{Schilling2}] 
  Let $\lambda^\circ=(\lambda^{\circ_1}, \ldots, \lambda^{\circ_k})$
  be a $k$-tuple of single row partitions, $\mu$ be a partition of
  weight $\vert \lambda^\circ\vert$ and $\delta$ be the partition such
  that $\delta_i=\vert \lambda^{\circ_i} \vert $. There exists a
  bijection $\Theta_k$ between ${\rm Tab}_{\lambda^\circ,\mu}$ and $RC(\mu,
  \delta)$ which is compatible with the cocharge and inversion
  statistic, i.e.
\begin{equation}
\forall\ T^{\circ} \in {\rm Tab}_{\lambda^{(k)}, \mu}\ , \ {\rm cc}\left 
(\Theta_k(T^\circ)\right ) = {\rm Inv}(T^\circ)\ .
\end{equation} 
\end{theorem}
We recall the steps of the algorithm permitting to compute
$\Theta_k(T^\circ)$. This algorithm is implemented in the package {\tt MuPAD-Combinat}.
\begin{algorithm}[\cite{Schilling2}]\label{algo} $ $ \\
  {\tt Input:} $T^\circ=(T^{\circ_1}, \ldots, T^{\circ_k})$ {\tt a
  $k$-tuple of Young tableaux of} ${\rm Tab}_{\lambda^\circ,\mu}$. \\\\
  {\tt Initialization: $\nu^\circ \longleftarrow$ a sequence of p empty partitions.} \\\\
  {\tt For i from k down to 1 do \\
    \indent For j from 1 to $l(\lambda^{\circ_i})$ do \\
    \indent \indent 1- For k from $T^{\circ_i}_j$ to p do \\
    \indent \indent \indent \indent Add a box in the j-th row in the partition $\nu^{\circ_k}$ \\
    \indent \indent \indent EndFor \\
    \indent \indent 2- Recompute all vacancy numbers,\\
    \indent \indent 3- Fill the new cells coming from step 1 \\
    \indent \indent \indent with the vacancy number of their row,\\
    \indent \indent 4- Remove a maximal number in the (j-1)-th row \\
    \indent \indent \indent of the partitions which have a new box from step 1 in the j-th row.\\
    \indent EndFor\\
    EndFor\\
    For a from 1 to p-1 do \\
    \indent For i from 1 to $\mu_1$ do \\
    \indent \indent Replace each number $\beta$ in the row $\nu^{\circ_a}_i$ by
 $\nu^{\circ_{a+1}}_i-\nu^{\circ_a}_i-\beta$\\
    \indent EndFor \\
    EndFor.  }
\end{algorithm}
%
\begin{proposition} 
  Let $k$ be a positive integer and $\lambda, \mu$ be two partitions
  such that $k\vert \mu \vert=\vert \lambda \vert$. Let
  $\lambda^{(k)}=(\lambda^{\circ_1}, \ldots, \lambda^{\circ_k})$ be the
  $k$-quotient of $\lambda$ and $\delta$ be the partition $(\vert
  \lambda^{\circ_1} \vert, \ldots, \vert \lambda^{\circ_k} \vert)$.
  Combining the bijection $\Theta_k$ and $\Psi_k$, we have
\begin{equation}
  \tilde G_{\lambda,\mu}^{(k)}(q) = \widetilde{I}_{\lambda^{(k)},\mu}(q) =
  \widetilde{S}_{\mu, \delta}(q) = \sum_{\nu^\circ\in C_{\mu, \delta}} 
  q^{\alpha(\nu^\circ)}
  \prod_{ {}_{1\le i \le \mu_1}^{1\le a \le p-1} }\left\lbrack
    {}_{\nu_i^{\circ_a}-\nu_{i+1}^{\circ_a},~
      \nu_i^{\circ_{a+1}}-\nu_{i}^{\circ_a}}^{~~~~~~~\nu_i^{\circ_{a+1}}-\nu_{i+1}^{\circ_a} }\right\rbrack\ \ .
\end{equation}
\end{proposition}
This proposition gives an explicit formula for the transition matrix
between Hall-Littlewood functions and monomials in terms of $q$-binomial coefficients.
\begin{corollary} 
  Let $\lambda^\circ$ be the $l(\lambda)$-quotient of the partition
  $l(\lambda)\lambda$ and $\delta$ be the partition defined by
  $\delta_i= \vert \lambda^{\circ_i} \vert$.
  Hall-Littlewood function can be expressed as
  \begin{equation}
    \widetilde{Q}^{'}_\lambda(X;q)=\sum_{\mu\dashv \vert \lambda \vert} \tilde{S}_{\mu, \delta}(q)\ m_\mu \ .
  \end{equation}
\end{corollary}
\subsubsection{A matricial recoding of the bijection}
We shall therefore propose a simpler but similar
algorithm for finding the shape of ribbon rigged configurations.
Let $p$ and $q$ be two integers and $\mathcal{M}_{p,q}$ the set of all
$(p\times q)$-matrices with integer coefficients.  We define the
operator $A_E$ on $\mathcal{M}_{p,q}$ by
\begin{equation}
\begin{array}{cccl}
A_E:   &        \mathcal{M}_{p,q}     & \longrightarrow & \mathcal{M}_{p,q} \\
       &    M=(m_{i,j})_{i,j}  & \longmapsto     & N = \left \lbrace \begin{array}{cccr}
                                                  n_{i,1} &=& m_{i,1} & \\
                                                  n_{i,j} &=& \sum_{1\le k \le j}{m_{i,k}} \ , \ \text{for} & 2\le j \le q\ .
                                                 \end{array} \right . 
\end{array}
\end{equation}
\noindent Let $T^\circ$ be a $k$-tuple of Young tableaux of shape $\lambda^\circ$, 
weight $\mu$ and diagonal vector $d$. Let
$\Theta_k(T^\circ)=(\nu_{T^\circ}, J_{T^\circ})$ be the ribbon rigged
configuration corresponding to $T^\circ$ by $\Theta_k$. We construct a
$\left (m_2\times l(\mu)\right)$-matrix $M^{T^\circ}$ using the
following rule
$$
M^{T^\circ}_{i,j}= \text{number of cells labelled $j$ in diagonal
  $d_{i+1}$}~.
$$
\begin{proposition}\label{shape}
  The $j$-th column of $A_E\left (M^{T^\circ}\right )$ is equal to the
  $j$-th partition of $\nu_{T^\circ}$.
\end{proposition}
{\it Proof:} Let $T^\circ$ be a $k$-tuple of single row tableaux and
$\Theta_k(T^\circ)=(\nu_{T^\circ}, J_{T^\circ})$. In the algorithm
\ref{algo}, we observe that boxes which appear in the $i$-th line of a
partition $\nu_{T^\circ}^{(j)}$ only come from elements of the $i$-th
diagonal of $T^\circ$ which are smaller than $j$. And by definition of
the operator $A_E$, the entry $(i,j)$ of the matrix $A_E(M^{T^\circ})$
corresponds to the number of cells less than $j$ in the $i$-th
diagonal. \hfill $\square$

\begin{example}
  Consider the following 3-tuple of single row tableaux $T^\circ$ 
$$\left ( \ 
\begin{array}{|c|c|}
\hline
1 & 4 \\
\hline
\end{array} \quad , \quad 
\begin{array}{|c|c|}
\hline
1 & 2 \\
\hline
\end{array} \quad , \quad
\begin{array}{|c|c|c|c|}
\hline
1 & 2 & 3 & 3 \\
\hline
\end{array}\ \right )\ .$$
In this case, the matrices $M^{T^\circ}$ and $A_E(M^{T^\circ})$ are 
\begin{center}
$$ M^{T^\circ}=\left ( \begin{array}{cccc} 3 & 0 & 0 & 0 \\
                                    0 & 2 & 0 & 1 \\
                                    0 & 0 & 1 & 0 \\
                                    0 & 0 & 1 & 0  
                 \end{array}\right ) \quad \text{and} \quad
 A_E(M^{T^\circ})=\left ( \begin{array}{cccc}  3 & 3 & 3 & 3 \\
                                     0 & 2 & 2 & 3 \\
                                     0 & 0 & 1 & 1 \\
                                     0 & 0 & 1 & 1
                 \end{array}\right ) \ .
$$
\end{center}
The shape of the rigged configuration $\Theta_k(T^\circ)$ is
\begin{center}
\includegraphics*[width=8cm]{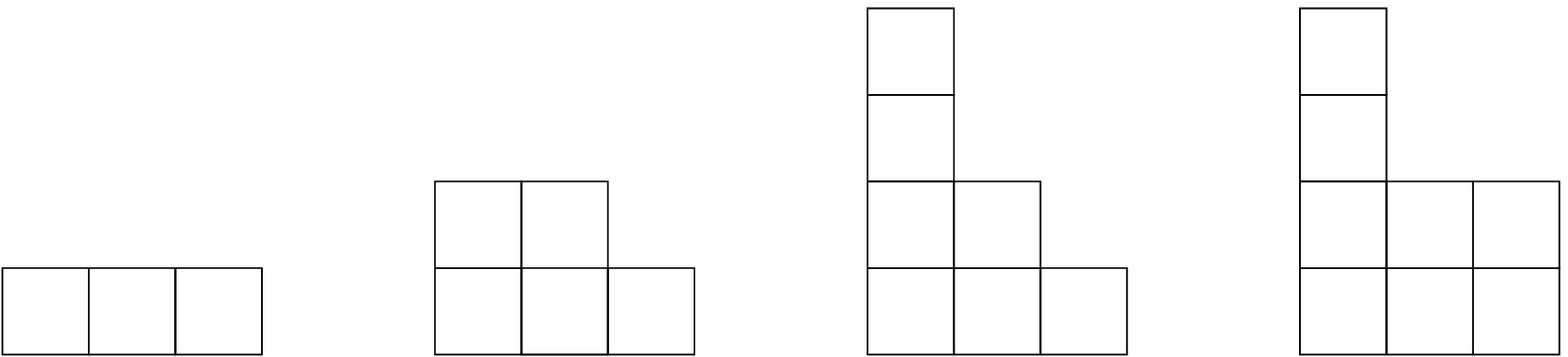}\ ,
\end{center}
and can be read from $A_E(M^{T^\circ})$.
\end{example}
\subsubsection{Diagonal classes on ribbon rigged configurations}\label{preuveP4}
We give an explicit formula for inversions polynomials restricted to diagonal
classes. This formula permits to prove the Proposition \ref{root}.

\begin{proposition}\label{propdiago}
  Let $\lambda^\circ$ be a sequence of single row partitions
  $(\lambda^{\circ_1}, \ldots, \lambda^{\circ_k})$, $\mu$ be a partition of
  weight $\sum_i\vert \lambda^{\circ_i} \vert$ and $\delta$ be the partition
  $(\vert \lambda^{\circ_1}\vert, \ldots, \vert\lambda^{\circ_k} \vert)$. For
  each diagonal vector in $\Delta_{\lambda^\circ, \mu}$, there
  exists a unique $(\mu, \delta)$-configuration $\nu^\circ$ in $C_{\mu, \delta}$
  such that
  $$ 
  \Theta_k(D_{\lambda^\circ, \mu}(d)) = RC_{\mu, \delta}(\nu^\circ)\ .
  $$
  The explicit expression for the inversion polynomial restricted to
  the diagonal class $D_{\lambda^\circ, \mu }(d)$ is
\begin{equation}
\widetilde{I}_{\lambda^\circ,\mu}(q;d)=\widetilde{S}_{\mu,
  \delta}(q;\nu^\circ)=q^{\alpha(\nu^\circ)} \prod_{ {}_{1\le i \le \mu_1}^{1\le a
    \le p-1} }\left\lbrack {}_{\nu_i^{\circ_a}-\nu_{i+1}^{\circ_a},~
    \nu_i^{\circ_{a+1}}-\nu_{i}^{\circ_a}}^{~~~~~~~\nu_i^{\circ_{a+1}}-\nu_{i+1}^{\circ_a}
  } \right\rbrack\ .
\end{equation}
\end{proposition}
\textit{Proof:} Let $d$ be a diagonal vector and $T^\circ$, $T^{'\circ}$ two
elements in the diagonal class $D_{\lambda^\circ,\mu}(d)$.  These two
$k$-tuples of tableaux differ only by a permutation of cells which are in a
same diagonal $d_i$. By construction, this property implies
$M^{T^\circ}=M^{T^{'\circ}}$ and by Proposition \ref{shape}, the ribbon
rigged configuration $\Theta_k(T^\circ)$ has the same shape $\nu^\circ$ than
$\Theta_k(T^{'\circ})$. Hence, since $\Theta_k$ is a bijection,
$D_{\lambda^\circ,\mu}(d)$ is embedded into $RC_{\mu, \delta}(\nu^\circ)$. \\
Conversely, let $T^\circ$ and $T^{'\circ}$ be two $k$-tuples of tableaux in
${\rm Tab}_{\lambda^\circ,\mu}$ which are not in the same diagonal class.
This implies that $M^{T^\circ}\not = M^{T^{'\circ}}$ and the shape of their
corresponding ribbon rigged configurations are not the same. Finally, we
conclude that
$$\Theta_k(D_{\lambda^\circ,\mu,}(d)) = RC_{\mu, \delta}(\nu)\ .$$
The expression of inversion polynomials of diagonal classes in terms of
$q$-supernomial coefficients follows immediately from the invariance of the
statistics under $\Theta_k$ and (\ref{FermionicDiago}).\hfill $\square$
\\\\

\begin{corollary}
  Let $n$ be a positive integer, $\lambda^\circ$ be the $k$-tuple of row
  partitions $((n),\ldots, (n))$ and $\mu$ a partition of weight $nk$
  satisfying the condition $\mu=(ks_1,\ldots, ks_p)$ for some positive
  integers $s_1,\ldots, s_p$ such that $s_1+\ldots+s_p=n$.  Let
  $D_{\lambda^\circ,\mu}(d)$ be a diagonal class with only one element
  $T^\circ$ and $\Theta_k(T^\circ)=(\nu_{T^\circ}, J_{T^\circ})$ its
  corresponding ribbon rigged configuration. The $i$-th partition of the shape
  $\nu_{T^\circ}$ is the rectangular shape $(k^{s_1+\ldots + s_i})$.
\end{corollary}
\textit{Proof:} By Proposition \ref{one}, since $T^\circ$ is alone in its
diagonal class, $T^\circ$ is a $k$-tuple of tableaux of shape $\lambda^\circ$
with the same values at the same positions of each single row tableau. The
corresponding matrix $M^{T^\circ}$ is
$$ M^{T^\circ}=\left ( \begin{array}{cccccc}  
 k      &     0  & \ldots & \ldots &\ldots     & 0       \\
\vdots  & \vdots &        &        &           & \vdots  \\
 k      &     0  & \ldots & \dots  &\ldots     & 0       \\
 0      &     k  &   0    & \ldots &\ldots     & 0       \\
 \vdots & \vdots & \vdots &        &           & \vdots  \\
 0      &     k  &   0    & \ldots & \ldots    & 0       \\
 \vdots &        &        &        &           & \vdots  \\
 0      &\ldots  & \ldots & \ldots &  0        & k       \\
 \vdots &        &        &        &  \vdots   & \vdots  \\
 0      & \ldots & \dots  & \ldots &  0        & k  
\end{array}\right )\ ,
$$
where $k$ occurs $s_i$ times in the $i$-th column. This implies
that the shape is given by the matrix 
$$ A_E(M^{T^\circ})=\left ( \begin{array}{cccccc}  
 k      &     k  & \ldots & \ldots &\ldots     & k       \\
\vdots  & \vdots &        &        &           & \vdots  \\
 k      &     k  & \ldots & \dots  &\ldots     & k       \\
 0      &     k  &   k    & \ldots &\ldots     & k       \\
 \vdots & \vdots & \vdots &        &           & \vdots  \\
 0      &     k  &   k    & \ldots & \ldots    & k       \\
 \vdots &        &        &        &           & \vdots  \\
 0      & \ldots & \ldots & \ldots &    0      & k       \\
 \vdots &        &        &        & \vdots    & \vdots  \\
 0      & \ldots & \ldots & \ldots &    0      & k  
 \end{array}\right )\ . 
$$
Then, the $i$-th partition in the shape $\nu_{T^\circ}$ is the
rectangular partition $(k^{s_1+\ldots + s_i})$.  \hfill $\square$
\\\\
Now, we are able to give a combinatorial proof of Proposition
\ref{root}. Let $T^\circ$ be a $k$-tuple of tableaux of shape
$((n),\ldots, (n))$ in a diagonal class $D_{\lambda^\circ,\mu}(d)$
with strictly more than one element.  Write
$\nu^\circ=(\nu^{\circ_1},\ldots,\nu^{\circ_p})$ the shape of the corresponding
ribbon rigged configuration $\Theta_k(T^\circ)$. By Proposition
\ref{propdiago}, this shape is the same for all $k$-tuples of tableaux
in $D_{\lambda^\circ,\mu}(d)$. Let $h$ be the last position such that the
$h$-th diagonal $d_h$ has at least two different elements. Then, the
$(h+1)$-th partition in $\nu_{T^\circ}$ is a rectangular partition of
width $k$ and height $s\le r$. Since the last part of
$\nu^{\circ_h}=(\nu^{\circ_h}_1,\ldots,\nu^{\circ_h}_l)$ is equal to $a$ with $a <
h$, the following coefficient appears in the inversion polynomial
$\widetilde{I}_{\lambda^\circ, \mu}(q;d)$
$$
\left \lbrack{}_{\nu^{\circ_h}_l-\nu^{\circ_h}_{l+1},~~
    \nu^{\circ_{h+1}}_l-\nu^{\circ_h}_{l}}^{~~~~~~~\nu^{\circ_{h+1}}_l-\nu^{\circ_h}_{l+1}}
\right \rbrack = \left \lbrack{}_{a-0,~k-a}^{~~~~k} \right \rbrack~.
$$
By definition, all $k$-th primitive roots of unity $\zeta$ are
roots of the $q$-binomial coefficient $\left \lbrack{}_{a,~k-a}^{~~~k}
\right \rbrack$. Finally
$$
\widetilde{I}_{\lambda^\circ, \mu}(\zeta; d) = 0\ .
$$ 
\hfill $\square$

\section{Conclusion}
The main tool of our combinatorial proof is the explicit expression of
inversion polynomials on diagonal classes in terms of $q$-binomial
coefficients. In this approach, we have used ribbon rigged configurations and
fermionic formulas given in \cite{Schilling2} which only exist in the case of
products of row partitions or column partitions. The unrestricted rigged
configurations (\cite{Schilling3}) are an other kind of configurations defined
in the case of products of rectangles. Unfortunately, in the special case of
product of rows, the number of shapes and the number of diagonal classes are
not the same.  Hence, the corresponding fermionic formula $M$ cannot be used.
We can mentionned \cite{Schilling4} for a survey of the zoology of rigged
configurations.
 \vspace{0.5cm}


\end{document}